Cite as: Gao, H., Evans, T., & Fergusson, A. (2025). Student Explanation Strategies in Tertiary Mathematics and Statistics Education: A Scoping Review. In S. Cook, B. Katz, & M.-R. D. (Eds.), Proceedings of the 27th Annual Conference on Research in Undergraduate Mathematics Education (pp. 258-268). SIGMAA on RUME. Alexandria, VA. https://doi.org/10.17608/k6.auckland.29925140

# Student Explanation Strategies in Postsecondary Mathematics and Statistics Education: A Scoping Review

Huixin Gao
University of Auckland

Tanya Evans
University of Auckland

Anna Fergusson
University of Auckland

*This scoping review examines the use of student explanation strategies in postsecondary mathematics and statistics education. We analyzed 46 peer-reviewed articles published between 2014 and 2024, categorizing student explanations into three main types: self-explanation, peer explanation and explanation to fictitious others. The review synthesizes the theoretical underpinnings of these strategies, drawing on the retrieval practice hypothesis, generative learning hypothesis, and social presence hypothesis. Our findings indicate that while self-explanation and explaining to fictitious others foster individual cognitive processes enhancing generative thinking, peer explanation have the potential to combine these benefits with collaborative learning. However, explanation to fictitious others have the potential to mitigate some of the negative impacts that may occur in peer explanation, such as more knowledgeable students dominating peer discussions. The efficacy of the methods varies based on implementation, duration, and context. This scoping review contributes to the growing body of literature on generative learning strategies in postsecondary education and provides insights for optimizing the integration of student explanation techniques within mathematics and statistics.*

## Introduction

In mathematics and statistics education, explanations serve as a critical tool for conveying complex concepts and fostering deep understanding (Rittle-Johnson & Loehr, 2017). These explanations can be broadly categorized into two types: teacher explanations and student explanations.

Teacher explanations traditionally dominate both face-to-face and online video-based instruction in postsecondary education, often criticized for fostering student passivity. In response, recent research has increasingly explored strategies to reduce the reliance on teacher-centered methods, particularly following Freeman et al.'s (2014) meta-analysis, which found that active learning significantly boosts student performance in Science, Technology, Engineering, Mathematics (STEM) fields. From a cognitive perspective, meaningful learning and long-term knowledge retention are shown to require active engagement in the process of knowledge generation and construction, rather than passive consumption (Fiorella & Mayer, 2014; Freeman et al., 2014; Rittle-Johnson et al., 2017)

Student explanations, on the other hand, involve learners articulating their understanding of concepts, processes, or problem-solving steps, either to themselves or others. Hodds et al. (2014) demonstrated that prompting undergraduate students to elucidate mathematical concepts enhances their focus on analyzing and explicating the learning content, thereby facilitating improved learning outcomes. In general, there is a consensus that student explanations play a crucial role in shaping understanding and engagement (Bisra et al., 2018; Fiorella & Mayer, 2014; Rittle-Johnson et al., 2017).

Despite the recognized importance of student explanations, current research still exhibits several gaps in integrating their theoretical bases. Additionally, the specifics of implementation processes of student explanations in mathematics and statistics remain relatively vague. This investigation aims to contribute to the growing body of literature on generative learning strategies in postsecondary educational contexts, by categorizing the literature about student explanation in mathematics and statistics education at the tertiary level.

## Theoretical grounding

The cognitive processes underlying student explanation can be understood through the lens of Human Cognitive Architecture (Atkinson & Shiffrin, 1968). This model posits that information processing occurs in stages: from sensory memory (visual, auditory, and tactile inputs) to working memory (where information is manipulated) and finally to long-term memory (where knowledge is stored and retrieved). These processes are fundamental to understanding how student explanations facilitate learning in mathematics and statistics education.

The retrieval practice hypothesis, proposed by Roediger & Karpicke (2006), aligns with the role of long-term memory in this cognitive model (Karpicke & Roediger, 2008). It posits that the act of explanation enhances learning through active retrieval of previously acquired knowledge. In mathematics and statistics education, this can be applied by encouraging students to explain complex formulas (Lyle et al., 2020), statistical concepts (Hun Lim et al., 2015), or problem-solving procedures without referring to notes or textbooks (May, 2022). This process strengthens memory traces and improves long-term memory, as demonstrated by Lachner et al. (2020, 2021).

The generative learning hypothesis (Lachner et al., 2018) relates to both working memory and long-term memory processes. It suggests that student explanation involves complex cognitive processes beyond mere retrieval. In higher-level mathematics and statistics courses, this can be implemented by having students create explanations that connect different mathematical concepts or apply statistical principles to novel problems. This approach emphasizes active knowledge construction and the development of meaningful connections between concepts, engaging both working memory for immediate processing and long-term memory for knowledge integration (Fiorella & Mayer, 2014).

The social presence hypothesis (Short et al., 1976) also involves working memory and long-term memory processes, but through an interpersonal dimension. In mathematics and statistics education, this can be applied through peer teaching activities or group problem-solving sessions. This theory suggests that anticipating an audience's knowledge leads to continuous adjustment of explanatory content, enhancing learning (Mills et al., 2020). Empirical evidence shows that learners engaged in student-to-student explanation reported higher levels of social presence (Nasir et al., 2023), motivation (Morris et al., 2023), and increased mental effort (Alegre et al., 2019), indicating enhanced engagement of working memory and potential for long-term memory consolidation.

These three hypotheses collectively offer a comprehensive framework for understanding the efficacy of student explanation in mathematics and statistics education, highlighting the interplay between sensory input, cognitive processing, and social factors in enhancing learning outcomes.

## Methods

The approach for the scoping review is underpinned by Arksey & O'Malley's (2005) five-stage framework, which adopts a rigorous process of transparency, enabling replication of the search strategy and increasing the reliability of the study findings. The five stages of the framework are: (1) identifying the initial research questions; (2) identifying relevant studies; (3) study selection; (4) charting the data; and (5) collating, summarizing and reporting the results. These stages were utilized in this review of the student explanation literature.

**Identifying the initial research questions**

The focus of our review was to examine the specifics of the use of student explanation in mathematics and statistics postsecondary education, thus we proposed the following research question:

*How can literature on student explanation strategies in postsecondary mathematics and statistics education be categorized to reveal current research patterns and trends?*

**Identifying relevant studies**

All publications included in this review were published in English and indexed in Eric, Scopus, or Web of Science databases. Additionally, we utilized websites (i.e., Google Scholar and Research Rabbit https://www.researchrabbit.ai/) and conducted citation searches to ensure comprehensive coverage. The linked descriptive key search terms that were developed to guide the search are outlined in Table 1. Research prior to the past decade, when educational reforms have incorporated principles of active learning into higher education curricula, is unlikely to reflect these specific teaching techniques in higher education. Therefore, the time period 2014 to 2024 was considered appropriate. The search strategy's exclusive focus on the term 'explanation' may have inadvertently excluded relevant studies that incorporated explanatory activities as intervention components but utilized different terminology in their documentation. A full list of inclusion and exclusion criteria is outlined in Table 2.

*Table 1. Search string*

| Structured combination of keywords, phrases, and Boolean operators |
| --- |
| (“Higher Education” OR “Postsecondary” OR “Tertiary Education” OR “Undergraduate” OR “College”) AND (“Mathematics” OR “Statistics” OR “Science education” OR “Engineering science” OR “STEM”) AND (“Learning by explaining” OR “Explaining-based learning” OR “Explanation-driven learning” OR "Self-explanation" OR "Self-explaining" OR “Learning by teaching” OR “Learning-by-teaching” OR “Team-based-explanation” OR “Peer tutoring” OR “Peer instruction” OR “Peer assist” OR “Peer-led” OR “Explanation to fictitious others” OR “ Non-interactive learning”) |

*Table 2. Inclusion and exclusion criteria*

| Criterion | Inclusion | Exclusion |
| --- | --- | --- |
| Time period | 2014 to 2024 | Studies outside these dates |
| Language | English | Non-English studies |
| Type of article | Original research, published in a peer review journal | Articles that were not peer-reviewed or published elsewhere |
| Literature focus | Articles where the overwhelming theme relates to student explanation | Articles that made a passing or token reference to student explanation |
| Population and sample | Undergraduate or college students enrolled in a recognized program of study | All other students not enrolled in undergraduate study |

**Study selection and data charting**

The process of study selection followed the Preferred Reporting of Items of Systematic Reviews and Meta-Analyses (PRISMA) protocol (Moher et al., 2009). We implemented a search strategy across Eric, Scopus, and Web of Science databases, which yielded 759 initial results. After removing duplicates and screening titles and abstracts, 102 publications progressed to the full-text evaluation stage. Based on our predefined criteria, we excluded 56 publications. Ultimately, 46 studies were included in the final analysis (the full list is available at figshare: https://doi.org/10.17608/k6.auckland.28202357). Summaries were developed of each study related to the author, year, location of study, study design, study methods, and sample size, and a brief comment on the limitations and recommendations of the individual selected study.

## Findings

Based on the analysis and synthesis of the results of selected literature, the student explanation strategies can be categorized into three broad types: self-explanation, peer explanation and explanation to fictitious others. Of the 46 studies reviewed, 15 focus on self-explanation (12 on mathematics, 3 on statistics), 26 discuss peer explanation (9 on STEM, 14 on mathematics, and 3 on statistics), while 5 mainly focus on explanation to fictitious others (3 on STEM, 2 on mathematics).

**Self-explanation, peer explanation, and explanation to fictitious others**

Self-Explanation (SE) is a cognitive process wherein individuals elucidate to themselves the meaning of new information, often integrating it with existing cognitive schema (Bichler et al., 2022). This strategy facilitates learners' active participation in knowledge construction (Hall & Vance, 2010). SE are generated for the learner, not intended to teach the content to other people (Rittle-Johnson et al., 2017). Specifically, these methods can be broadly categorized into three main forms: written (e.g., documenting explanations on paper or in digital text boxes) (Salle, 2020), oral (e.g., recording verbal explanations via audio devices or self-narration) (Hodds et al., 2014), and gestural (e.g., conveying information through body language) (Salle, 2020).

The 15 SE studies we reviewed primarily explore the efficacy of SE. These studies examined various aspects, including promoting learning outcomes (e.g., Atisso et al., 2023; Rittle-Johnson & Loehr, 2017; Ryan & Koppenhofer, 2022), self-explanation for proof comprehension (Hodds et al., 2014; Maarif et al., 2020; Wirth et al., 2024), the impact of self-explanation quality on learning (Bichler et al., 2022), types of self-explanation and their impact on learning (Salle, 2020), and influence on cognition (Hadfield, 2021; Miller-Cotto & Auxter, 2021).

As postsecondary education has evolved, peer tutoring and cooperative learning approaches have been increasingly integrated with explanatory techniques in university settings (Naughton et al., 2021; Pollock et al., 2023). This study employs the term "Peer Explanation" (PE) to refer specifically to the explanatory process within peer instruction at undergraduate and postgraduate levels. PE focuses on explanatory interactions among students, emphasizing the impact of collaborative dynamics on the quality and effectiveness of explanations. At the same time, PE seems to include the process of SE, since students convey their own understanding to their peers. Unlike the independent cognitive process of SE, PE requires real-time bi-directional or multi-directional communication, demanding immediate question-and-answer interactions between the explainer and the audience, where listeners can raise questions, request clarification, or seek additional information (Maxwell & Wiles, 2022). PE typically involves a formalized process where students engage in mutual explanation within learning groups of at least three participants. In classroom settings, PE typically operates in pairs or fixed groups (3-5 members), with clearly defined roles (explainer and listeners) and structured time allocations and rotation mechanisms.

The 26 PE studies we reviewed examine various aspects, including impact on learning motivation (Mayerhofer et al., 2023), improvement of academic achievement (e.g., Biju, 2019; Corbo & Sasaki, 2021; Han & Thiel, 2022; Kishbaugh et al., 2018; Maxwell & Wiles, 2022; Morales et al., 2014; Nasir et al., 2023), effect of mixed-gender activities (Dumitru et al., 2022; Powell et al., 2023), cognitive development (Kustina & Meredith, 2023), the relationship with self-work during lectures (Leitner & Gabel, 2024), and the development of self-efficacy (Mills et al., 2020; Morris et al., 2023; Naughton et al., 2021). Advantages and Disadvantages of Online PE (Campbell, 2019).

With technological advances transforming postsecondary education, innovative methods for capturing student explanations through audio and video recordings have emerged. These developments have given rise to a new learning activity known as Explanation to Fictitious Others (EFO) (Lachner et al., 2018). EFO refers to students generating explanations and then explaining them to non-present, fictitious others via video, audio or text. While EFO also involves student-generated explanations, it differs from SE in that the explanatory content must be recorded through media tools (video, audio, or text), making it available for sharing and self-refinement (Hoogerheide et al., 2016). Generally, it includes three stages (Fiorella & Mayer, 2014): first, students prepare their explanations; second, they have to explain the subject matter to potential others either orally or in writing; third, students answer subsequent questions in interactive settings or receive feedback from an instructor. EFO triggers distinct cognitive processes that are essential for students' knowledge acquisition (Lachner et al., 2018).

The 5 EFO studies we reviewed mainly focus on comparing the effects of explaining through writing and explaining via video (Hoogerheide et al., 2016; Lachner et al., 2018), the impact of EFO on learning achievements and students' self-perceptions (Fiorella & Mayer, 2014; Hoogerheide et al., 2014), the effects of the timing of EFO implementation on students (Lachner et al., 2020), and the impact of EFO on students' problem-solving performance (Hoogerheide et al., 2019). EFO typically yields the most significant results when students record videos explaining problems (Hoogerheide et al., 2016, 2019).

**Theoretical underpinnings of SE, PE and EFO**

The Retrieval Practice Hypothesis and Generative Learning Hypothesis provide theoretical support for the efficacy of SE and EFO in enhancing students' long-term memory through retrieval and more advanced generative learning mechanisms such as organization and integration (Lachner et al., 2021). This is particularly evident in the learning performance of students with different levels of prior knowledge (Hodds et al., 2014; Maarif et al., 2020; Ryan & Koppenhofer, 2022; Leppink et al., 2012). Low prior knowledge students perform worse on challenging questions due to insufficient knowledge reserves, which hinders their working memory. Therefore, PE can offer a useful knowledge link to address this disparity. Because of differences in students' knowledge, PE could be beneficial to further balance their learning. The Social Presence Hypothesis provides additional theoretical support for PE. The students' interactions and discussions with each other are able to some extent promote their cognitive thinking and build their understanding (Mills et al., 2020; Morris et al., 2023). This collaborative aspect of PE could help bridge the gap between students with different levels of prior knowledge, facilitating a more balanced learning environment where students can benefit from each other's strengths and compensate for individual weaknesses.

However, PE may also have limitations. Students with lower levels of understanding might unintentionally reinforce each other's misconceptions, leading to the spread of incorrect information (Kishbaugh et al., 2018). Additionally, the quality of explanations can vary significantly, as not all students possess the ability to clearly articulate complex concepts (Bichler et al., 2022). Without proper guidance and oversight from the instructor, peer-led discussions may sometimes fail to address critical gaps in understanding, potentially hindering overall learning outcomes (Smith et al., 2009). Furthermore, students who are more confident or knowledgeable may dominate the discussion, limiting the contributions of others and reducing the collaborative benefits intended by peer explanation (Mills et al., 2020).

Thus, EFO, by eliminating real-time peer interactions, mitigates some of the potential negative effects that the presence of peers can have in peer explanations. Although the target audience for the explanation is non-present fictitious others, it still involves the stages of preparation and explanation. In the preparing stages, EFO encourages learners to engage in a multi-faceted cognitive process that seems to differ to SE. With the expectation to teach others, this process involves selecting salient information to be conveyed, organizing it into a coherent mental representation, and integrating it with existing knowledge structures. This aligns with the teaching expectancy effect, as described by Fiorella and Mayer (2013), which posits that anticipation of teaching enhances one's own learning.

The explanation stage leverages the power of perspective-taking. Learners, when constructing explanations for their imagined audience, are compelled to consider diverse viewpoints. This cognitive exercise not only enhances the quality of their explanations but also reinforces their own understanding. As Hoogerheide et al. (2016) suggest, the belief that one's explanations may influence others can serve as a potent catalyst for self-directed learning. Moreover, the concept of non-present others in EFO serves a dual purpose. It heightens learners' awareness of their potential audience, thereby increasing their level of arousal and engagement (Hoogerheide et al., 2019). This heightened state of cognitive alertness may contribute to more effective information processing and retention.

## Conclusions

This scoping literature review synthesizes the theoretical underpinnings of student explanation and categorizes student explanations into three main types: SE PE, and EFO. While SE and EFO foster individual cognitive processes that support effective generative thinking, PE combines these benefits with collaborative learning. However, EFO have the potential to mitigate some of the negative impacts that may occur in PE, such as more knowledgeable students dominating peer discussions. The retrieval practice, generative learning, and social presence hypotheses underpin the effectiveness of these approaches. However, the varying efficacy of SE, PE and EFO warrants further exploration, as current research reports mixed outcomes depending on factors such as implementation, duration, and context. Future research should focus on consolidating findings through systematic or integrative reviews to better understand and compare the impact of SE, PE and EFO on student learning outcomes.